\documentclass[12pt]{article}
\usepackage{amsmath,amssymb,stmaryrd}
\usepackage[dvips]{graphicx}

\begin{document}

\date{April 20, 2007}
\title{\Large\bf Smoothed Wigner transforms in the numerical simulation of 
semiclassical (high-frequency) wave propagation. }
\author{Agissilaos G. Athanassoulis$^1$ \\[2mm]
 $^1$The Program in Applied and Computational Mathematics,\\
 Princeton University \\ Email: aathanas@math.princeton.edu\\[2mm]
} \maketitle

{\footnotesize \noindent {\bf Abstract.}  The numerical simulation of wave propagation in semiclassical (high-frequency) problems is well known to pose a formidable challenge. In this work, a new phase-space approach for the numerical simulation of semiclassical wave propagation, making use of the smoothed Wigner Transform (SWT), is proposed.  There are numerous works which use the Wigner Transform (WT) in the study of a variety of wave propagation problems including high-frequency limits for linear, nonlinear and/or random waves. The WT however is well known to present significant difficulties in the formulation of numerical schemes.  Working with concrete examples for the semiclassical linear Schr\"{o}dinger equation it is seen that the SWT approach is indeed significantly {\em faster} (in a well-defined sense) to work with than the WT and than full numerical solutions of the original equation in the semiclassical regime. Comparisons with exact and numerical solutions are used to keep track of numerical errors. \\

\noindent
{\bf Keywords.} Smoothed Wigner Transform, semiclassical wave propagation, phase-space, particle methods, Schr\"{o}dinger equation \\
{\bf AMS (MOS) subject classification:} 65T99, 65M99, 35S10 } \vskip.2in

\section{Introduction}

\noindent Phase space methods for wave propagation, such as these presented here, are in the junction between time-frequency analysis \cite{Co1,Co2,Co3,Fl1,Fl2,Hla,Jaf,Jan,Mal}, microlocal analysis \cite{Ben,Ge1,Ge2,Ji1,Spa,Zha}, PDEs \cite{Bal,Bao,Eng,Fan,Fil} and physics \cite{Gro,Hal,Moy,Tat,Wig} (at least). Distributions from Cohen's class (i.e. Wigner Transforms convolved with various smoothing kernels) have been long well known and used in signal processing contexts, as semi-empirical ways to `clean up' the time-frequency structure of various chirps, e.g. \cite{Fl1}. Microlocal measures have evolved in the last twenty years, in close connection with PDEs, in fact with a very similar job in mind, e.g. \cite{Ge1,Ge2}. The need to extend powerful spectral-density-based ideas from stationary to nonstationary stochastic problems gave birth to various nonstationary spectral densities \cite{ Bal,Fan,Fl2,Hal,Ryz}. These three lines of work all converge to the Wigner Transform (WT), and explain the recent explosion in the literature for WT methods for wave propagation (see e.g.\cite{Ath,Bal,Bao,Ben,Eng,Fan,Fil,Ge2,Hal,Ji1,Ji2,Ji3,Lio,Mar,Mau,Ryz,Spa,Zha}). One soon discovers however that, despite the analytic and asymptotic power of these methods, there has been much less work in computing with WTs (for numerical analysis with WTs see \cite{Arn,Bal,Bao}). This is not an accident, and many people complain that it is not at all practical to compute with WTs, see e.g. \cite{Ben,Jaf}. Computations with Wigner measures are feasible, but this is a different thing, and we discuss it below. In this work we propose a way to go around the WTs `bad' features, and exploit its `good' features in the numerical solution of semiclassical wave propagation problems. To do that, we work with a particular distribution from Cohen's class, which we will call the Smoothed Wigner Transform (SWT).

This paper is a continuation of \cite{Ath}. The two papers are complementary in the following way: in \cite{Ath} the results corresponding to individual numerical experiments were presented in some detail. The objective there was to present direct evidense that `the SWT approach works'. However, as is discussed in section \ref{sec4_2}, the efficiency of the method can only be compared to the more conventional approaches when working with a population of numerical experiments for different values of the semiclassical parameter $\varepsilon$. In this work we focus on this synoptic question, and present clear evidence that the SWT approach is indeed more efficient -- in a well defined sense -- than more conventional approaches. \\

Historically, the WT
\[
W^\varepsilon \left[ {f,g} \right]\left( {x,k} \right)=\int\limits_{y\in 
\mathbb{R}^n} {e^{-2\pi iky}f\left( {x+\frac{\varepsilon y}{2}} \right)\bar 
{g}\left( {x-\frac{\varepsilon y}{2}} \right)dy} 
\]
was introduced for the study of the semiclassical limit of quantum mechanics \cite{Wig}, i.e. the high-frequency Schr\"{o}dinger equation. Consider the system of linear equations

\begin{equation}
\label{eq1}
\displaystyle
\begin{array}{l}
\varepsilon \frac{\partial }{\partial t}u^\varepsilon \left( {x,t} 
\right)+L^\varepsilon \left( {x,\varepsilon \partial _x } 
\right)u^\varepsilon \left( {x,t} \right)=0,
\\
\left( {u_i^\varepsilon \left( {x,0} \right)} \right)_{i=1,...,d} =\left( 
{u_{i;0}^\varepsilon \left( x \right)} \right)_{i=1,...,d}
\end{array}
\end{equation}

\noindent for the complex-valued $d$-component wavefield $u^\varepsilon \left(x,t\right)=\left(u_i ^\varepsilon \left(x,t\right) \right)_{i=1,...,d} $, where $x\in \mathbb{R}^n$  and the $d\times d$ matrix-valued Weyl symbol $L^\varepsilon \left( {x,k} \right)$ of the spatial pseudodifferential operator $L^\varepsilon \left( {x,\varepsilon \partial _x } \right)$ is anti-Hermitian, $L^\varepsilon \left( {x,k} \right)=- \overline { L^\varepsilon \left( {x,k}  \right)}^T $. This is a natural way to quantify the requirement that equation (\ref{eq1}) is a wave equation; indeed many interesting linear wave equations have anti-Hermitian Weyl symbols. The use of the Weyl pseudodifferential calculus with WTs is natural, see e.g. \cite{Gro,Poo}. The small parameter $\varepsilon $ scales the typical wavelength of the wavefunctions; the length scale at which the coefficients vary (or equivalently one over the slopes of the coefficients), as well as the typical propagation distance we're interested in are significantly larger than $\varepsilon $. Many semiclassical problems in quantum mechanics and long distance propagation problems in (classical) continuum mechanics can be cast in the form (\ref{eq1}) \cite{Ge2,Ryz,Whi}. For the remainder of this paper when we refer to a high-frequency wave propagation problem we will mean a problem that can be cast in the form (\ref{eq1}) with anti-Hermitian Weyl symbol.

There can be (at least) two problems associated with the IVP (\ref{eq1}): The \textit{high-frequency-limit problem}: recover the limits of certain observables of the wavefield, such as the $L^2$-norm density,
\[
N^0\left( {x,t} \right):=\mathop {\lim }\limits_{\varepsilon \to 0} \left| {\sum\limits_{i=1}^d {u_i^\varepsilon \left( {x,t} \right)} } \right|^2
\]
 and the energy density 
\[
E^0\left( {x,t} \right):=\mathop {\lim }\limits_{\varepsilon \to 0} \,\,\left\langle {L^\varepsilon \left( {x,\varepsilon \partial _x } \right)u^\varepsilon \left( {x,t} \right),u^\varepsilon \left( {x,t} \right)} \right\rangle _{\mathbb{R}^d}.
\]
 On the other hand, the \textit{high-frequency problem} consists in solving for the observables at a given nonzero $\varepsilon \ll 1$; naturally, simplifications making use of the fact that $\varepsilon $ is small can be made, but the idea is not to lose all the $\varepsilon $-dependent 
information. The difference is nontrivial; for example it is often the case that limit-observables have singularities (become infinite) while for any $\varepsilon >0$ the same observables stay finite for all spacetime points.

WT methods for high-frequency wave propagation problems consist in studying 
the Wigner matrix $W^\varepsilon \left({x,k,t} \right) =  \left[ {W^\varepsilon \left[ {u_i^\varepsilon ,u_j^\varepsilon } \right]\left( {x,k,t} \right)} \right]_{i,j=1,...,d} $ instead of the wavefields $\left( {u_i^\varepsilon \left( {x,t} \right)} \right)_{i=1,...,d} $ themselves. It can be said that this is a good idea for two reasons: first of all, due to `time-frequency density' character of the WT wave-like functions are in principle `nicely' represented in phase space. Moreover, the operator for the evolution in phase space is 
(asymptotically) simpler, and becomes a Hamiltonian flow in many cases. In 
fact this is a general theme in WT-based methods; there are `two 
transforms': one of the `data' $ {u^\varepsilon \left( {x,t} \right)}  \mapsto  {W^\varepsilon \left( {x,k,t} \right)} $, and one of the `dynamics' under which the data evolve.

A formidable tool for high-frequency-limit problems consists in Wigner 
measures (WM), \[W^0\left( {x,k} \right):=\mathop {w- \lim }\limits_{\varepsilon 
\to 0} \left[ {W^\varepsilon \left[ {u_i^\varepsilon ,u_j^\varepsilon } 
\right]\left( {x,k} \right)} \right]_{i,j=1,...,d} \] \cite{Ge1,Lio}. Pretty much all bilinear observables of the wavefields can be expressed (in the limit) 
in terms of the WM. Recently computational methods for the high-frequency-limit wave propagation based on the WM have been developed \cite{Ji1,Ji2,Ji3}.

Although often discussed in the context of the WT, the question of smoothing phase-space dynamic models and its implications has only rarely been actually studied in a deterministic context - see however \cite{Co1,Mur}.

The remainder of this paper is structured as follows: in section 2 we review the WM approach for the limit problem; SWTs and the phase-space reformulations of problem (\ref{eq1}) in terms of them are reviewed in section 3. Numerical results on the SWT methods for the high-frequency problem are presented in section 4; conclusions and further work are discussed in section 5.

\section{The WM approach for the high-frequency-limit problem}

\noindent For an excellent (defining) overview of the WM approach one should refer to \cite{Ge2} and for a comparison with WKB-based asymptotic techniques to \cite{Spa}. WMs are a refinement of microlocal concepts related to microlocal defect measures, H-measures etc \cite{Ge1,Ge2,Lio}. They provide a description for `asymptotic signals'; properly speaking, an asymptotic signal is a family of functions $\left\{ {u^\varepsilon \left( x \right)} \right\}_{\varepsilon >0} $ that in some sense share the same time-frequency structure. A prototype asymptotic signal is given by the WKB ansatz 
\[
\left\{ {u^\varepsilon \left( x \right)=A\left( x \right)e^{\frac{2\pi i}{\varepsilon }S\left( x \right)}} \right\}_{\varepsilon >0} , 
\]
and the corresponding WM is 
\[
W^0 \left( {x,k} \right) =  \left|{A\left( x \right)} \right|^2\delta \left( {k-\nabla _x S\left( x \right)} \right) ,
\]
 unfolding amplitude and `local wavenumber' (often referred to as `instantaneous frequency' in time-frequency analysis) information in phase space. Naturally, a unique `local wavenumber' doesn't always exist, as e.g. for multi-WKB signals $u^\varepsilon \left( x \right)=\sum\limits_{n=1}^N {A_n \left( x \right)e^{\frac{2\pi i}{\varepsilon }S_n \left( x \right)}} $. In fact an aspect of the WKB ansatz is that the number of `branches' $N$ at each point $\left( {x,t} \right)$ is a `hidden' parameter. This becomes very important when seeking a solution of (\ref{eq1}) in the form $u^\varepsilon \left( {x,t} \right)=\sum\limits_{n=1}^N {A_n \left( {x,t} \right)e^{\frac{2\pi 
i}{\varepsilon }S_n \left( {x,t} \right)}} $, because at different points 
of spacetime $\left( {x,t} \right)$ there are different unknown numbers of 
(nonzero) branches $N=N\left( {x,t} \right)$. The fact that WMs treat this (i.e. the `number of different local wavenumbers') in a more subtle way is an important advantage in the study of high-frequency limits \cite{Spa}.

For multicomponent wavefields the WM is a positive-definite-Hermitian-matrix-valued measure on phase space \cite{Ge2}. To 
link matrix-valued WMs with physical quantities in general one must work in the context of a wave problem like (\ref{eq1}). In that case 
\[
N^0\left( 
x \right)=\int\limits_{k\in \mathbb{R}^n} {tr\left( {W^0 \left( {x,k} 
\right)} \right)dk}
\]
is the total (i.e. taking into account all components 
of the wavefield) limit norm density and 
\[
E^0\left( x 
\right)=\int\limits_{k\in \mathbb{R}^n} {tr\left( {L^0\left( {x,k} 
\right)W^0 \left( {x,k} \right)} \right)dk}
\]
the total limit energy density [17,35]. Observe that the WT allows the construction of phase space 
densities, `elaborating' on the more traditional physical space 
densities, i.e. 
\[
N^0\left( {x,k} \right)=tr\left( {W^0 \left( {x,k} 
\right)} \right), \,\,\,
 E^0\left( {x,k} \right)=tr\left( {L^0\left( {x,k} 
\right)W^0 \left( {x,k} \right)} \right).
\]

The derivation of exact equations for the evolution in time of the WT $W^\varepsilon \left( {x,k,t} \right)$ corresponding to a problem of the form (\ref{eq1}) is a well studied topic \cite{Ge2,Lio,Mar,Mau,Moy,Ryz,Wig,Zha}. From them, asymptotic equations governing the evolution of WMs $W^0 \left( {x,k,t} \right)$ can be obtained. In simple cases they consist in decoupled Liouville equations in phase space for scalar phase-space measures \cite{Ge2,Lio}. Therefore, WT-based numerical methods for the high-frequency-limit problem deal with the numerical solution of the Liouville equation in phase space when the initial condition is an appropriate measure. Moment and level set methods have been formulated and used in this connection \cite{Eng,Ji1,Ji2,Ji3}.

A natural question arises here: can't these methods be used for the small 
$\varepsilon >0$ problem? So far the answer seems to be negative. The 
$\varepsilon \to 0$ weak limit `kills' a lot of fine scale features. If they 
are not killed, but included in the numerical problem, we get an overly 
heavy problem, that typically is slower to solve than a full numerical 
solution of (\ref{eq1}), thus beating its purpose. The need to put the extra data under control clearly emerges; this is possible in terms of the SWT.

\section{The SWT approach for the high-frequency problem}

\subsection{The WT vs. the SWT}

\noindent The Wigner Transform is defined as a sesquilinear mapping, 
\[
 W^\varepsilon :f\left( x \right),g\left( x \right)\mapsto W^\varepsilon \left[ {f,g} \right]\left( {x,k} \right). 
\]
When $f=g$ it is called the Wigner distribution of $f$, and denoted as $W^\varepsilon \left[ f \right]\left( {x,k} \right)$. The Wigner distribution of a wavefunction $f\left( x \right)$ is a good way to realize a joint breakdown of the wavefunction's $L^2$-norm over space $x$ and wavenumber $k$ - with the caution that it takes on negative values as well.  The books \cite{Fl1,Fol} are two very important sources on the WT and its properties, \cite{Fl1} emphasizing the signal processing point of view, and \cite{Fol} the relation between the WT and pseudodifferential operators.

Cohen's class of distributions \cite{Co2,Co3,Fl1,Jaf,Jan} is defined as the class of all sesquilinear transforms of the form 
\[
C\left[ {f,g} \right]\left( {x,k} \right)= \int\limits_{y,u,z} {f\left( {u+\frac{y}{2}} \right)\bar {g}\left( {u-\frac{y}{2}} \right)e^{-2\pi i\left[ {ky+zx-zu} \right]}\phi _C \left( {z,y} \right)dydudz}
\]
for any distribution $\phi _C \left( {z,y} \right)$ (restricting $f,g$ to test functions). An equivalent definition is the class of transforms that results from convolving the WT with a distributional kernel 
\[
K_C \left( {x,k} \right)=\mathcal{F}_{\left( {z,y} \right)\to \left( {x,k} 
\right)}^{-1} \left[ {\phi _C \left( {z,y} \right)} \right]
\]
\[
C\left[ {f,g} \right]\left( {x,k} \right)=\int\limits_{{x}',{k}'\in 
\mathbb{R}^n} {K_C \left( {x-{x}',k-{k}'} \right)\,\,W\left[ {f,g} 
\right]\left( {{x}',{k}'} \right)d{x}'d{k}'}.
\]

Among them there are many attractive alternatives to the WT. When the kernel 
is chosen to be the Wigner distribution of some function $h$, $K_C \left( {x,k} 
\right)=W\left[ h \right]\left( {x,k} \right)$, and if $f=g$, we get the 
spectrogram of $f$ with window $h$. Smoothing WTs with an appropriate kernel 
tames the interference terms, but doesn't necessarily kill them completely; 
there is a balance between smoothing enough and not smoothing too much. In 
this connection we define the scaled SWT as the sesquilinear transform

\begin{equation}
\label{eq2}
\begin{array}{c}
\widetilde{W}^{\sigma _x ,\sigma _k ;\varepsilon }\left[ {f,g} \right]\left( 
{x,k} \right) = \\
=\frac{2}{\varepsilon \sigma _x \sigma _k 
}\int\limits_{{x}',{k}'\in \mathbb{R}^n} {e^{-\frac{2\pi \left| {x-{x}'} 
\right|^2}{\varepsilon \sigma _x^2 }-\frac{2\pi \left| {k-{k}'} \right|^2}{\varepsilon \sigma _k^2 }}W^\varepsilon \left[ {f,g} 
\right]\left( {{x}',{k}'} \right)d{x}'d{k}'}. 
\end{array}
\end{equation} \\

This is a (scaled) WT convolved with a tensor-product Gaussian with space-domain variance proportional to $\varepsilon \sigma _x^2 $ and wavenumber-domain variance proportional to $\varepsilon \sigma _k^2 $. The scaling is selected with the problem (\ref{eq1}) and the semiclassical regime in mind; see Figures 1,2 for intuition on the scale of smoothing. The SWT provides  coarse-scale\footnote{Different authors favor different terms, e.g. slow-scale, coarse-scale, macroscopic, homogenized etc. Each of these terms stands for a description which contains information for the underlying object, understood to be sufficient in terms of certain qualitative criteria, but which is not complete -- usually it is essentially less than complete, thus leading to compression as well. } description of the time-frequency structure of the underlying signal / wavefield; its resolution, i.e. how much information is effectively discarded, can be tuned through the smoothing parameters.

\begin{figure}
\centering
\includegraphics[width=110mm,height=50mm]{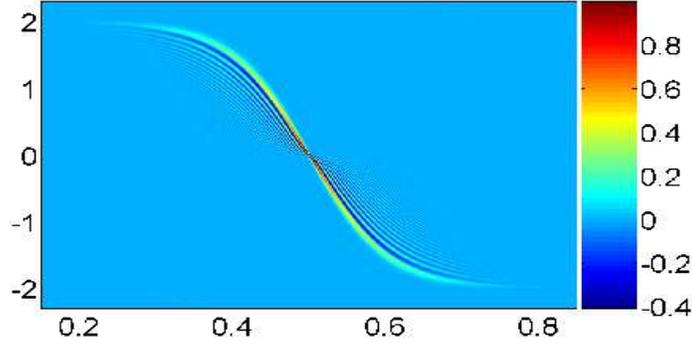}
\caption{Wigner Transform of the `tanh chirp', $f^{\varepsilon} (x) = A(x) e^{\frac{2\pi i}{\varepsilon} S(x)}$, $A(x)=e^{-25(x-0.5)^2}$, $S(x)=\frac{-1}{5} log \left( e^{10(x-0.5)}+e^{-10(x-0.5)} \right) $ 
{\footnotesize (Remark: Image quality has to do with -- the very understanable -- arXiv size limits, and does not represent the quality of the output; see e.g. \cite{Ath}}).
}
\end{figure}

\begin{figure}
\centering
\includegraphics[width=110mm,height=50mm]{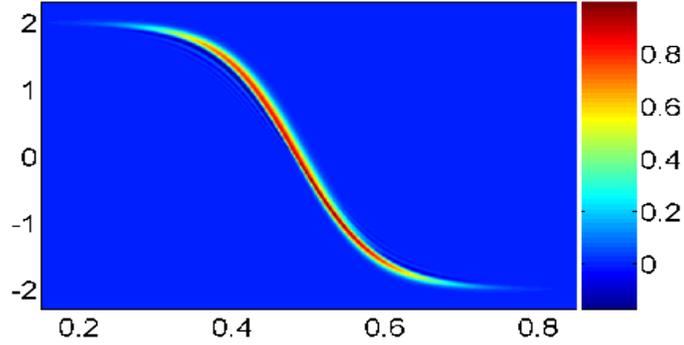}
\caption{Smoothed Wigner Transform of the `tanh chirp', correspondig to the scaling $\sigma_x,\sigma_k=O(1)$ 
{\footnotesize (Remark: Image quality has to do with -- the very understanable -- arXiv size limits, and does not represent the quality of the output; see e.g. \cite{Ath}}).
}
\end{figure}

\subsection{Equations for SWTs}

In \cite{Ath} (an early version of which can be found in the arXiv) we presented the derivation of exact equations for the SWT of a 
wavefunction. That result essentially stated that, given the multicomponent 
problem (\ref{eq1}), its smoothed Wigner matrix satisfies
\begin{equation}
\label{eq3}
\varepsilon \frac{\partial }{\partial t}\widetilde{W}+A+B^T=0,
\end{equation}
where

\[
A=L^\varepsilon \left( {x-\frac{\varepsilon }{4\pi i}\partial _k 
-\frac{\varepsilon \sigma _x^2 }{4\pi }\partial _x ,k+\frac{\varepsilon 
}{4\pi i}\partial _x -\frac{\varepsilon \sigma _k^2 }{4\pi }\partial _k } 
\right)\widetilde{W},
\]
\[
B=\overline{L^\varepsilon} \left( {x+\frac{\varepsilon }{4\pi i}\partial _k 
-\frac{\varepsilon \sigma _x^2 }{4\pi }\partial _x ,k-\frac{\varepsilon 
}{4\pi i}\partial _x -\frac{\varepsilon \sigma _k^2 }{4\pi }\partial _k } 
\right)\widetilde{W}^T, 
\]
\noindent with initial data 
\[
 \widetilde{W}\left( {x,k,0} \right)=\left[ 
{\widetilde{W}^{\sigma _x ,\sigma _k ;\varepsilon }\left[ 
{u_{i;0}^\varepsilon ,u_{j;0}^\varepsilon } \right]\left( {x,k} \right)} 
\right]_{i,j=1,...,d}. 
\] 

We have used the abbreviation $\widetilde{W}=\left[ {\widetilde{W}^{\sigma _x ,\sigma _k ;\varepsilon }\left[ {u_i^\varepsilon ,u_j^\varepsilon } \right]\left( {x,k,t} \right)} \right]_{i,j=1,...,d} $, and, as we discussed earlier, the Weyl PDO calculus. For concreteness we 
will now focus on the (scalar) linear Schr\"{o}dinger equation on $\mathbb{R} $, i.e.

\begin{equation}
\label{eq4}
\begin{array}{c}
\varepsilon \frac{\partial }{\partial t}u^\varepsilon \left( {x,t} 
\right)-i\frac{\varepsilon ^2}{2}\Delta u^\varepsilon \left( {x,t} 
\right)+iV\left( x \right)u^\varepsilon \left( {x,t} \right)=0,\\
u^\varepsilon \left( {x,0} \right)=u_0^\varepsilon \left( x \right)
\end{array}
\end{equation}

\noindent In that case (\ref{eq3}) becomes
\begin{equation}
\label{eq5}
\begin{array}{c}
 \frac{\partial }{\partial t}\widetilde{W}\left( {x,k,t} \right)+\left( 
{2\pi k\frac{\partial }{\partial x}+\frac{\varepsilon \sigma _k^2 
}{2}\frac{\partial ^2}{\partial x\partial k}} \right)\widetilde{W}\left( 
{x,k,t} \right)+ \\ 
+\frac{2}{\varepsilon} Re\left( 
{iV\left( {x-\frac{\varepsilon }{4\pi i}\partial _k -\frac{\varepsilon 
\sigma _x^2 }{4\pi }\partial _x } \right)\widetilde{W}\left( {x,k,t} 
\right)} \right)=0, \\ \\
\widetilde{W}\left( {x,k,t} \right)=\widetilde{W}^{\sigma _x ,\sigma _k 
;\varepsilon }\left[ {u_0^\varepsilon } \right]\left( {x,k,t} \right). \\
 \end{array}
\end{equation}

This is the exact equation, not an asymptotic approximation. Observe that the smoothing parameters $\sigma _x ,\sigma _k $ are now part of the data of the phase-space equation. `Little' smoothing, i.e. a regularized WT, corresponds to $\sigma _x ,\sigma _k \leqslant 1$ in the scaling (\ref{eq2}). In that case the asymptotic phase space equation is

\begin{equation}
\label{eq6}
\begin{array}{l}
 \frac{\partial }{\partial t}\widetilde{W}\left( {x,k,t} \right)+\left( 
{2\pi k\frac{\partial }{\partial x}-\frac{{V}'\left( x \right)}{2\pi 
}\frac{\partial }{\partial k}} \right)\widetilde{W}\left( {x,k,t} \right)+ 
\\ 
+\varepsilon \left( {{V}''\left( x 
\right)\frac{\sigma _x^2 }{8\pi ^2}+\frac{\sigma _k^2 }{2}} 
\right)\frac{\partial ^2}{\partial x\partial k}\widetilde{W}\left( {x,k,t} 
\right)=O\left( {\varepsilon ^2} \right). \\ 
 \end{array}
\end{equation}

\noindent In particular, the leading order part is a Liouville equation, and is the same as the corresponding leading-order equation for the WT.

The numerical solution of these equations for the SWT are expected to provide a way for a very efficient coarse-scale solver for wave propagation problems. Indeed, in the more standard approach of solving the original wave equation -- in this case equation (\ref{eq4}) -- we need to keep track of the wavefunction `exactly'. That is, even if a coarse-scale amplitude is all we are looking for, in order to be able to predict how the coarse-scale amplitude will evolve, we need as well to keep track of the fine-scale features of the wavefunction -- up to numerical errors of course. The SWT approach on the other hand, enables us to formulate (exactly or asymptoticly) a numerical problem for a coarse-scale description of the wavefield. Thus, we have to keep track of significantly less information, and ultimately we are able to carry out the computations significantly faster.

\section{Numerical results}
\subsection{The numerical scheme}

\noindent The subject of this paper is the formulation of a numerical method suitable for the simulation of semiclassical problems in terms of the SWT, and the investigation of its efficiency as compared to more conventional methods for the same problem. The test cases we work with is the IVP (\ref{eq4}), up to time $t_{max}$, independent of $\varepsilon$, with initial data scaled with $\varepsilon$, for different initial conditions and potentials. The model we discretize is the Liouville equation for the SWT with $\sigma _x ,\sigma _k \leqslant 1$, i.e. the leading order part of (\ref{eq6}). The conclusions of this investigation, based on the numerical results of sections 4.2 and 4.3 are discussed in section 5. An outline of the algorithm follows.

The WT is computed on a Cartesian grid in phase space with the FFT, with complexity $O\left( {N^2\log N} \right)$ \cite{Mal}. The complexity for the computation of the SWT is $O\left( {L^2N^2\log N} \right)$ where $L$ is the number of sampling points needed for the smoothing kernel. When $\sigma _x ,\sigma _k =O\left( 1 \right)$, $L$ is of $O\left( 1 \right)$. Adaptive computation of the SWT that doesn't spend much time on regions of phase space with no energy is also possible \cite{Ozd}, and might be essential for two- and higher dimensional problems.

The Liouville equation is solved numerically with the use of particles, i.e. the numerical implementation of the method of characteristics (related, but not identical, to \cite{ Arn,Fer,Rav}). An initial population of particles is created, so that the SWT can be interpolated up to an error tolerance from its values on them. The trajectory of each particle is computed according to Hamilton's ODEs with a Runge-Kutta solver; the value of the density on each particle remains unchanged in time. The solution at each moment in time is constructed by interpolating the density from its values on the particles.

\subsection{SWT computations vs. solutions of the original problem (\ref{eq4})}
\label{sec4_2}

The results of the SWT method are compared with exact as well as numerical solutions of the corresponding Schr\"{o}dinger equation (\ref{eq4}). Being able to compute faster with the SWT than by numerical methods for the original problem (\ref{eq4}) (or (\ref{eq1}) in general), is essential to our point of view. So far, the method we present here appears to be significantly faster than WT methods and numerical solutions of (\ref{eq4}); however first we have to make precise what exactly we mean by that.

It is well known that as the semiclassical parameter $\varepsilon$ becomes smaller, the numerical solution of (\ref{eq4}) eventually becomes intractable. Saying that a numerical method is {\em slower} or {\em faster} than another in the $\varepsilon \ll 1$ regime we mean the computation time for the `same' problem grows slower or faster in $\varepsilon$ for the one method than for the other. Hence we're looking for the {\em semiclassical complexity} $T=T\left( \varepsilon \right)$, that is the asymptotic dependence of the total computation time on the semiclassical parameter. In the same way, we're intersted in the rate of growth of the degrees of freedom needed for each numerical method $D=D(\varepsilon)$. That way we can compare numerical methods that are implemented in very different ways, and still get a well defined, meaningful result.

As was shown in \cite{Bao}, meaningful numerical solutions by means of finite differences methods for (\ref{eq4}) can be obtained with two different meshing strategies: a finer meshing strategy is needed for the correct $L^2$ approximation of the wavefunction, while a coarser meshing strategy can lead to correct slow-scale observables (but wrong fast-scale information due to numerical dispersion). As we discussed earlier, the SWT method can only recover slow-scale information. Independent numerical experiments, carried out by Kostas Politis (NTUA), lead to findings consistent with those of \cite{Bao}: the semiclassical complexity of standard finite differences methods appears to be $O \left( \varepsilon ^{-3} \right)$ or more, depending on the meshing strategy. 

A wavelet method developed by K. Politis shows a rate of growth of computation time around $O(N^{2.5})$, where $N$ is the number of points (i.e. translated versions of the finest-scale wavelets) used; $N$ is proportional to $\varepsilon^{-1}$. Thus wavelet methods seem very promising for semiclassical computations as well. Wavelet based solutions also appear to be better behaved than finite differences in terms of numerical dispersion, and can be accelerated with thresholding.

In the numerical experiments carried out so far, the SWT approach successfully captures the slow-scale observables of the wavefunction with semiclassical complexity no more than $O \left( -\varepsilon ^{-2} log \left( \varepsilon \right) \right) $; see Figures 3, 4. It is in that sense that our preliminary investigation shows that the SWT method is {\em faster} than full solutions of the problem (\ref{eq4}). (In all examples considered the SWT was significantly faster than the wavelet method, although the rates of growth of computation time in $\varepsilon^{-1}$ are relatively close for the two methods). A more thorough study of this question of semiclassical complexity, including more theoretical aspects as well as different numerical methods (e.g. time-splitting spectral methods, which seem to be particularly well adapted to semiclassical computations) is in order, and in fact in progress, in collaboration with other groups as well.

\begin{figure}
\centering
\includegraphics[width=110mm,height=70mm]{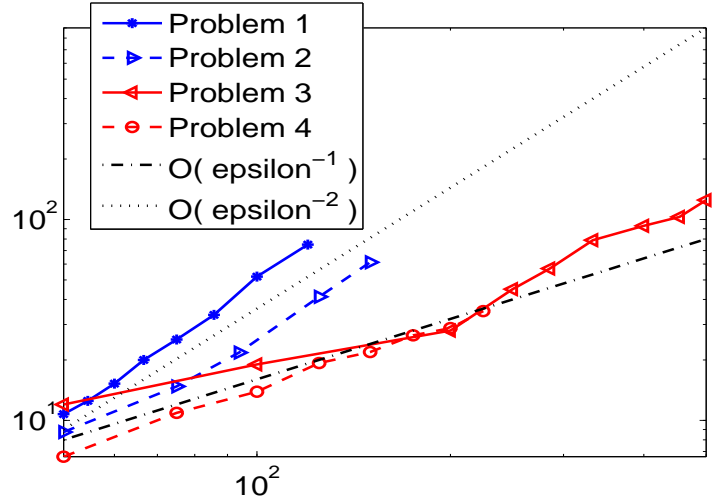}
\caption{Total computation time against $\varepsilon^{-1}$. All problems are special cases of (4); Problem 1 corresponds to
$u_0^{\varepsilon} (x) = A(x) e^{\frac{2\pi i}{\varepsilon} S(x)} $, 
$ A(x)=\frac{1}{4} \left(  tanh\left( 6.87(x+2.42) \right) +1  \right)$
$ \left(  tanh\left( 6.87(2.42-x) \right) +1  \right)$,
$S(x)=-\frac{x^4}{4}-x^2+2x$; $V(x)=0$;
Problem 2 corresponds to $u_0^{\varepsilon} (x) = A(x) e^{\frac{2\pi i}{\varepsilon} S(x)} $, $A(x)$ as above, $S(x)=-\frac{x^4}{4}+2x$; $V(x)=x$;
Problem 3 corresponds to $u_0^{\varepsilon} (x) = e^{-\frac{1+7i}{10\varepsilon}x^2}
+e^{-\frac{0.2+3i}{10\varepsilon}x^2}
+e^{-\frac{0.9-8i}{10\varepsilon}x^2}$; $V(x)=x$;
Problem 4 corresponds to $u_0^{\varepsilon} (x) = e^{-(1+\frac{3i}{\varepsilon})x^2-2x-4}+e^{-(1+\frac{2i}{\varepsilon})x^2-x-1}
+e^{-(1+\frac{i}{\varepsilon})x^2-\frac{2}{3}x-\frac{4}{9}}$; $V(x)=0$. The SWT solutions are found to be valid slow-scale representations of full numerical solutions of Problems 1,2,3 (courtesy of K. Politis), and of the exact solutions of Problem 4. Quantifying in a satisfactory way how good a `slow-scale representation' is has proved somewhat subtle, and is in progress; see \cite{Ath} for qualitative comparisons. The time for the computation of the SWT of the initial condition, although asymptotically important, is very small in this examples.  Certain normalizations were necessary to fit all the data in the same axes; in any case the slopes are not affected.
}
\end{figure}

\begin{figure}
\centering
\includegraphics[width=110mm,height=70mm]{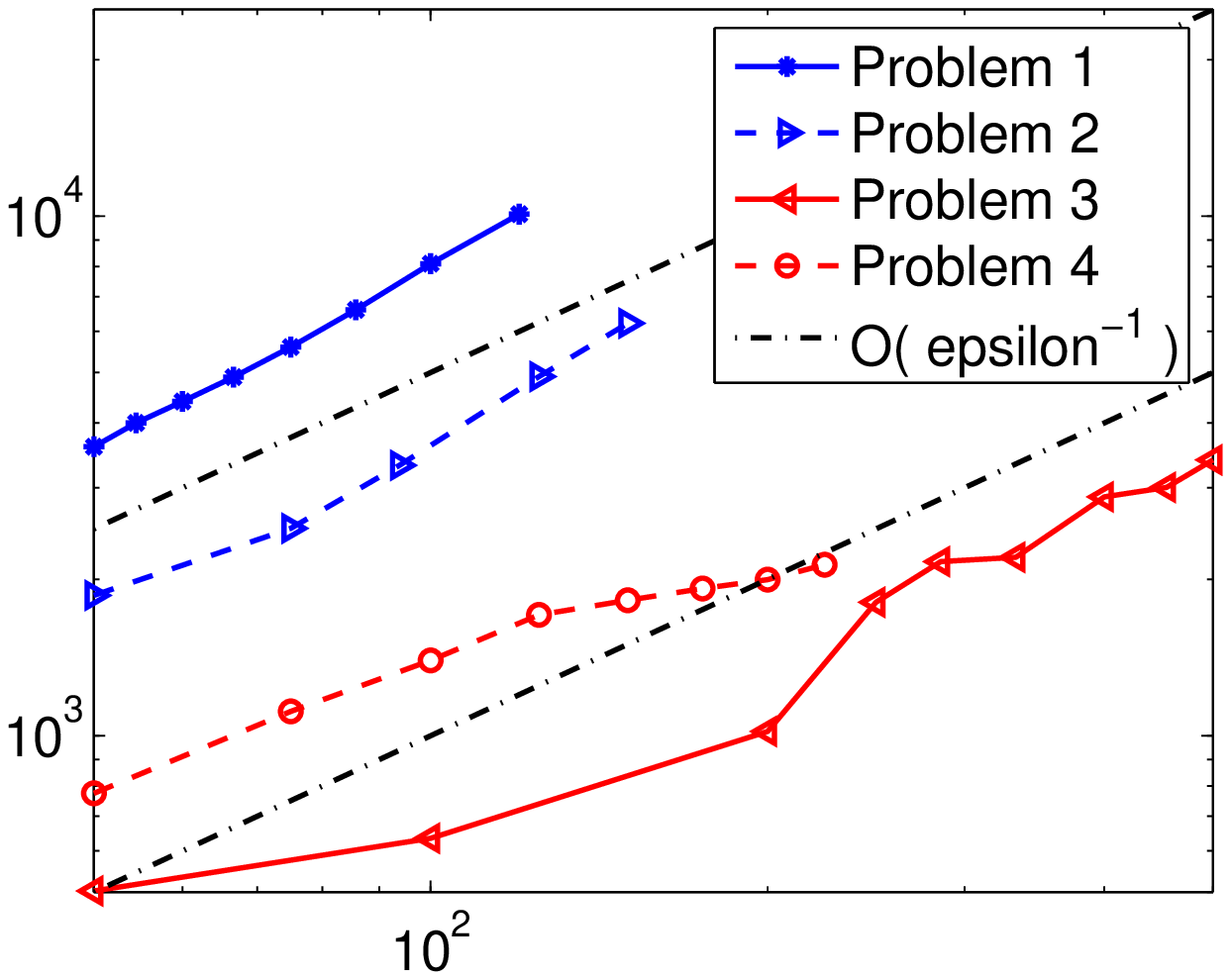}
\caption{Number of particles used against $\varepsilon^{-1}$.}
\end{figure}

\subsection{SWT vs. WT computations}

WT methods are slower than either full solutions of (\ref{eq4}) or SWT solutions. As compared with the SWT method we propose, the only difference would be another, more complicated initial condition. Apart from being represented with more points, the initial WT is qualitatively different than the SWT. It features many short waves with important amplitude, being a high-frequency 2-d wavefunction itself. 

Our numerical experiments indicate that WT methods are in any case one order of magnitude in $\varepsilon^{-1} $ slower, and sometimes even more. One particular phenomenon we observed in those numerical experiments is that, especially in problems with strong interference terms, in order to recover a meaningful approximation of the evolution of the WT many more points are required than in order to represent it initially. That is, if we use just enough particles to approximate the initial WT well, then often our numerical solution evolves to something very noisy; we need to use many more particles than that. This point is in contrast with the respective behavior for the SWT, and in fact we have no precise understanding at this point of how the number of particles for a given WT should be selected. (When this behavior kicks in, the complexity of the WT solution seems to be more than $O \left( \varepsilon ^{-3} \right)$).

\section{Conclusions -- Further work}
\noindent Semiclassical computations are well known to be a very demanding numerical problem. However, it seems that the more precise study of how different methods, in general solving for different things, behave as the semiclassical parameter goes to zero has started only recently. In this work we carry out a preliminary numerical investigation of a recently formulated SWT-based method for semiclassical simulations \cite{Ath}, and find that it is in general faster than full solutions of the original (physical space) PDE, based on numerical experiments, consistently with analysis of the more traditional methods \cite{Bao}. Also, that -- as expected -- SWT based computations are significantly faster and better behaved than WT computations.

As we mentioned earlier, equation (\ref{eq5}) is in fact a whole family of phase-space models, which can then be asymptotically treated to create a hierarchy of approximations for each phase-space model. In this paper we have only worked in a specific asymptotic regime, the semiclassical regime. One reason for this is that it allows direct comparisons with existing WT-based studies; also the history and evolution of phase-space methods draws heavily on semiclassical problems. However it should be clear that the principal reason for the SWT being `faster' is the fact that it is, and remains as time passes, a much simpler function than the WT. In other words {\em homogenization} here comes from the smoothing, and it is an essentially different thing than {\em asymptotics}.

These insights lead to the following conclusion: numerical methods corresponding to the exact equations should be formulated and compared to full solutions for different problems (e.g. `any' linear problem) and various smoothing choices; we can then look for practical, possibly problem specific, `slow-scale solvers'. Thus, further work at hand includes the treatment of different models included in (\ref{eq5}). Numerical methods for such models are needed; parallel implementation might be crucial for higher-dimensional problems.

\section{Acknowledgements}

The author would like to thank Profs. G. Papanicolaou, I. Daubechies, G. Makrakis, and L. Ryzhik for helpful discussions. Also, K. Politis for helpful discussions and for providing numerical results used for comparison. This work has been partially supported by NSF grant DMS-0530865.


\footnotesize

\begin{thebibliography}{99}

\bibitem{Arn}
A. Arnold, F. Nier, Numerical analysis of the deterministic particle method applied to the Wigner equation, {\em Math. Comp.} {\bf 58}, (1992) 645-669
\bibitem{Ath}
A. Athanassoulis, Exact equations for smoothed Wigner transforms and homogenization of wave propagation, {\em submitted for publication}
\bibitem{Bal}
G. Bal, L. Ryzhik, Time splitting for wave equations in random media, {\em Math. Model. Numer. Anal.} {\bf 38}, (2004) 961-987
\bibitem{Bao}
W. Bao, S. Jin, P. A. Markowich, On time-splitting spectral approximations for the Schr\"{o}dinger equation in the semiclassical regime, {\em J. Comput. Phys.}, {\bf 175}, (2002) 487-524.
\bibitem{Ben}
J. D. Benamou, F. Castella, T. Katsaounis, B. Perthame, High frequency limit of the Helmholtz equation, {\em Rev. Mat. Iberoamericana}, {\bf 18} (2002) 187-209
\bibitem{Co1}
L. Cohen, L. Galleani, Nonlinear transformation of differential equations into phase space, {\em EURASIP J. Appl. Signal Process.} {\bf 12} (2004) 1770-1777 
\bibitem{Co2}
L. Cohen, Time-frequency analysis, Prentice Hall PTR, 1995
\bibitem{Co3}
L. Cohen, Generalized phase space distributions, {\em J. Math. Phys.} {\bf 7} (1966) 781-786
\bibitem{Eng}
B. Engquist, O. Runborg, Computational high frequency wave propagation, {\em Acta Numer.} {\bf 12} (2003) 181-266 
\bibitem{Fan}
A. C. Fannjiang, S. Jin, G. Papanicolaou, High frequency behavior of the focusing nonlinear Schr\"{o}dinger equation with random inhomogeneities, {\em SIAM J. Appl. Math.} {\bf 63} (2003) 1328-1358 (eetron)
\bibitem{Fer}
E. Fergadakis, Numerical experiments with the particle method for the Wigner equation in high-frequency paraxial approximation, {\em unpublished}
\bibitem{Fil}
S. Filippas, G. N. Makrakis, Semiclassical Wigner function and geometrical optics, {\em Multiscale Model. Simul.} {\bf 1} (2003) 674-710 (eetron)
\bibitem{Fl1}
P. Flandrin, Time-frequency/time-scale analysis, Academic Press, San Diego CA, 1999
\bibitem{Fl2}
P. Flandrin, W. Martin The Wigner-Ville spectrum of nonstationary random signals, in: W. Mecklenbrauker, F. Hlawatsch (Eds), The Wigner Distribution, Elsevier, Amsterdam, 1997, pp. 211-267
\bibitem{Fol}
G. B. Folland, Harmonic analysis in phase space, Princeton University Press, Princeton NJ, 1989
\bibitem{Ge1}
P. Gerard, Microlocal defect measures, {\em Comm. Partial Differential Equations}, {\bf 16} (1991) 1761-1794 
\bibitem{Ge2}
P. Gerard, P. A. Markowich, N. J. Mauser, F. Poupaud, Homogenization limits and Wigner transforms, {\em Comm. Pure Appl. Math.} {\bf 50} (1997) 323-379 
\bibitem{Gro}
Grossmann, G. Loupias, E. M. Stein, An algebra of pseudodifferential operators and quantum mechanics in phase space, {\em Ann. Inst. Fourier (Grenoble)}, {\bf 18} (1968) 343-368
\bibitem{Hal}
B. Hall, M. Lisak, D. Anderson, R. Fedele, V. E. Semenov, Statistical theory for incoherent light propagation in nonlinear media, {\em Phys Rev E Stat. Nonlin. Soft Matter Physics}, {\bf 65} (2002)
\bibitem{Hla}
F. Hlawatsch, P. Flandrin, The interference structure of the Wigner distribution and related time-frequency signal representations, in: W. Mecklenbrauker, F. Hlawatsch (Eds), The Wigner Distribution, Elsevier, Amsterdam, 1997, pp. 59-133
\bibitem{Jaf}
S. Jaffard, Y. Meyer, R. D. Ryan, Wavelets, SIAM, Philadelphia PA, 2001
\bibitem{Jan}
J. E. M. Janssen, Positivity and spread of bilinear time-frequency distributions, in: W. Mecklenbrauker, F. Hlawatsch (Eds), The Wigner Distribution, Elsevier, Amsterdam, 1997, pp. 1-58
\bibitem{Ji1}
S. Jin, X. Li, Multi-phase computations of the semiclassical limit of the Schr\"{o}dinger equation and related problems: Whitham vs. Wigner, {\em Phys. D}, {\bf 182} (2003) 46-85 
\bibitem{Ji2}
S. Jin, H. Liu, S. Osher, R. Tsai, Computing multi-valued physical observables for the high frequency limit of symmetric hyperbolic systems, {\em J. Comput. Phys.} {\bf 210} (2005) 497-518
\bibitem{Ji3}
S. Jin, H. Liu, S. Osher, R. Tsai, Computing multivalued physical observables for the semiclassical limit of the Schr\"{o}dinger equation, {\em J. Comput. Phys.} {\bf 205} (2005) 222-241 
\bibitem{Lio}
P. Lions, T. Paul T, Sur les mesures de Wigner, {\em Rev. Mat. Iberoamericana}, {\bf 9} (1993) 553-618 
\bibitem{Mal}
S. Mallat, A wavelet tour of signal processing, Academic Press, 1999
\bibitem{Mar}
P. A. Markowich, N. J. Mauser, F. Poupaud, A Wigner-function approach to (semi)classical limits: Electrons in a periodic potential, {\em J. Math. Phys.} {\bf 35} (1994) 1066-1094 
\bibitem{Mau}
N. J. Mauser, (Semi)classical limits of Schr\"{o}dinger-Poisson systems via Wigner transforms, in: Journees Equations Aux Derivees Partielles, Ep. No. 12, Univ. Nantes, 2002
\bibitem{Moy}
J. E. Moyal, Quantum mechanics as a statistical theory, {\em Proc. Cambridge Philos. Soc.} {\bf 45} (1949) 99-124 
\bibitem{Mur}
C. Muratov, On the well-posedness of equations for smoothed phase space distribution functions and irreversibility in classical statistical mechanics, {\em J. Phys. A.} {\bf 34} (2001) 4641-51 
\bibitem{Ozd}
K. Ozdemir, O. Arikan, Fast computation of the ambiguity function and the Wigner distribution on arbitrary line segments, {\em IEEE Trans. Signal Process.} {\bf 49} (2001) 381-393 
\bibitem{Poo}
J. C. T. Pool, Mathematical aspects of the Weyl correspondence, {\em J. Math. Phys.} {\bf 7} (1966) 66-76
\bibitem{Rav}
P.-A. Raviart, An analysis of particle methods, in: Numerical Methods in Fluid Dynamics (Como, 1983); Lecture Notes in Math. 1127, Springer (1985) 243-324, 
\bibitem{Ryz}
L. Ryzhik, G. Papanicolaou, J. B. Keller, Transport equations for elastic and other waves in random media, {\em Wave Motion} {\bf 24} (1996) 327-70 
\bibitem{Spa}
C. Sparber, P. A. Markowich, N. J. Mauser, Wigner functions versus WKB-methods in multivalued geometrical optics, {\em Asymptot. Anal.} {\bf 33} (2003) 153-187
\bibitem{Tat}
V. I. Tatarskii, The Wigner representation of quantum mechanics, {\em Sov. Phys. Usp.} {\bf 26} (1983) 311-327
\bibitem{Whi}
G. B. Whitham, Linear and nonlinear waves, John Wiley, New York, 1999
\bibitem{Wig}
E. Wigner, On the quantum correction for thermodynamic equilibrium, {\em Phys. Rev.} {\bf 40} (1932) 749
(2006) 184-215 
\bibitem{Zha}
P. Zhang, Y. Zheng, N. J. Mauser, The limit from the Schr\"{o}dinger-Poisson to the Vlasov-Poisson equations with general data in one dimension, {\em Comm. Pure Appl. Math.} {\bf 55} (2002) 582-632

\end{thebibliography}
\end{document}